\documentstyle[art10]{article}

\title{Dimension of the Moduli Space of curves with an involution.}
\author{Luis Fuentes Garc\'{\i}a\thanks{Supported by an F.P.U.
fellowship of Spanish Government}
 \and{Manuel Pedreira P\'erez}
 }
\date{}

\newtheorem{teo}{Theorem}[section]

\newtheorem{prop}[teo]{Proposition}
\newtheorem{cor}[teo]{Corollary}
\newtheorem{lemma}[teo]{Lemma}

\newtheorem{rem}[teo]{Remark}

\font\euf=eufm10 at 12pt

\def\Curv{{\cal C}}
\def\CurvHyper{{\cal C}h}

\def\b{\mbox{\euf b}}
\def\bb{\mbox{{\euf b}}}

\def\R{{\cal B}}

\def\k{{\cal K}}
\def\K{{\cal K}}

\def\P{{\bf P}}

\newcommand\Te{{\cal O}}
\newcommand\C{C}

\newcommand\B{\beta}

\def\qed{\hspace{\fill}$\rule{2mm}{2mm}$}
\newcommand\lrw{\longrightarrow}

\begin{document}
\maketitle

{\footnotesize{\bf Authors' address:} Departamento de Algebra, Universidad de Santiago
de Compostela. $15706$ Santiago de Compostela. Galicia. Spain. e-mail: {\tt
pedreira@zmat.usc.es}; \\ {\tt luisfg@usc.es}\\
{\bf Abstract:}  Given a smooth curve $X$ of genus $g$ we compute de dimension of the family of curves
$C$ which have an involution over $X$. Moreover we distinguish when the curve $C$ is
hyperelliptic.
\\ {\bf Mathematics Subject Classifications (1991):} Primary, 14H10; secondary, 14H30, 14H37.\\ {\bf Key
Words:} Curves, involution.}

\vspace{0.1cm}

{\Large\bf Introduction.}  
Let ${\cal M}_{\pi}$ be the Moduli of smooth curves of genus $\pi$. It is well known that a generic
smooth curve $C\in {\cal M}_{\pi}$ with $\pi\geq 3$ does not have nontrivial automorphism. In particular
a generic smooth curve does not have involutions.

However, given an smooth curve $X$ of genus $g$ and a divisor $\b\in Div(X)$ verifying that $2(\b-\K)$ is
smooth we can construct a double cover $\gamma:C\lrw X$, where $C$ is a smooth curve of genus $\pi\geq
2g-1$. $C$ is not generic because it has an involution induced by $\gamma$. Moreover, we can characterize
the divisors $\b$ which provides double covers
$\gamma:C\lrw X$ with $C$ hyperelliptic (see \cite{fuentes2}). Thus if we denote by
$\Curv_{\pi}^{g}$ the family of curves of genus $\pi$ with an involution of genus
$g$ and by $\CurvHyper_{\pi}^{g}$ the family of hyperelliptic curves of genus $\pi$ with an involution of
genus
$g$, we have that:
$$
\CurvHyper_{\pi}^{g}\subset \Curv_{\pi}^{g}\subset {\cal M}_{\pi}.
$$
In this paper we compute the dimension of $\Curv_{\pi}^{g}$ and $\CurvHyper_{\pi}^{g}$. Let
$r=2(\pi-1)-4(g-1)$. The main results are:
\begin{enumerate}

\item If $r\geq 0$ then $dim(\Curv_{\pi}^{g})=2\pi-g-1$. Moreover, the dimension of the family of curves
with an involution over a hyperelliptic curve of genus $g$ is $2\pi-2g+1$.

\item If $r<0$ then $dim(\Curv_{\pi}^{g})=\emptyset$.

\item If $r=0,2,4$ then $dim(\CurvHyper_{\pi}^{g})=\pi$.

\item If $r<0$ or $r>4$ then $\CurvHyper_{\pi}^{g}=\emptyset$.
\end{enumerate}

Moreover, we proof that $\CurvHyper_{\pi}^{g}=\CurvHyper_{\pi}^{\pi-g}$. In fact we prove that given an
hyperelliptic curve $C$ of genus $\pi\geq 2$, an involution of  $C$ of genus $g$ composed with
the canonical involution provides an involution of genus $\pi-g$.

\bigskip

\section{Preliminaries.}

Let $C$, $X$ be two smooth curves of genus $\pi$ and $g$ respectively. Let $\gamma:C\lrw X$ be a double
cover. We know the following facts (see \cite{fuentes2}):
\begin{enumerate}

\item $\gamma_*\Te_C\cong \Te_X\oplus \Te_X(\K-\b)=S_{\bb}$ is a decomposable geometrically ruled
surface over the curve $X$. We call it canonical geometrically ruled surface. $\b$ is a nonspecial
divisor on $X$ verifying $2\b-2\k\sim \R$, where $\R$ is the branch divisor. If $X_0$ is the curve of
minimum self-intersection of $S_{\bb}$ and $X_1\sim X_0+(\b-\k)f$ then $C\sim 2X_1$. In particular $C\in
\langle 2X_0+\R f, 2X_1 \rangle$.

\item Conversely, let $\Te_X\oplus \Te_X(\K-\b)$ be a decomposable geometrically ruled surface over the
curve $X$, such that $2\b-2\k$ is a smooth divisor. Then the generic curve $C$ in the linear system
$|2X_1|$ is smooth and then we have a double cover $C\lrw X$. Moreover, if $C\in \langle 2X_0+\R f, 2X_1
\rangle$, with $\R\sim 2\b-2\k$ then the branch divisor of the cover is
$\R$.

\end{enumerate}

In order to compute the dimension of the curves with an involution over a fixed curve $X$, we have to
study when two curves in the linear system $|2X_1|$ are isomorphic. In this way we have the following
proposition:

\begin{prop}\label{unicidad}

Let $X$ be a smooth curve of genus $g$. Let $\b$ be a nonspecial divisor of degree $b\geq 2g-2$ defining
a canonical ruled surface $S_{\bb}$. Let $\R\sim 2\b-2\k$ be different points. Then, there is a unique
curve $\C\in |2X_1|$ up to isomorphism with a $2:1$ map $\gamma:C \lrw X$ whose ramification points
over $X$ are the points of $\R$.

\end{prop}
{\bf Proof:} We know that given a curve
$X_1\in |X_0+(\b-\k)f|$, a generic curve
$\C$ of the pencil $L=\langle 2X_0+\R f, 2X_1 \rangle$ is a curve with an involution $\gamma:\C\lrw X$
and ramification points at
$\R$ (see Lemma $1.9$ and Theorem $1.10$ in \cite{fuentes2}). This curve is invariant by the unique
involution of
$S_{\bb}$ that fixes the curves
$X_0$ and
$X_1$. In this way, $\C$ meets each generator in two points and these points are related by the
involution.

Let $\C$ and $\C'$ be two curves of the pencil $L$. Let $Pf$ be a generic generator. We can define an
automorphism of $S_{\bb}$ that fixes $X_0\cap Pf$, $X_1\cap Pf$ and takes a point of $\C\cap Pf$ into a
point of $\C'\cap Pf$. If we consider the restriction of this automorphism to the pencil $L$ we see
that it takes $\C$ into $\C'$ and the two curves are isomorphic.

Now, let $X_1$ and $X_1'$ two irreducible curves of the linear system $|X_0+(\b-\k)f|$. Since $S_{\bb}$
is a decomposable ruled surface we can define an automorphism of $S_{\bb}$ that takes $X_1$ into
$X_1'$. In this way the curves of the pencils $L=\langle 2X_0+\R f, 2X_1 \rangle$ and $L=\langle 2X_0+\R
f, 2X_1'
\rangle$ are isomorphic and our claim follows. \qed

\section{Computing the dimensions.}

\begin{prop}\label{dim1}

Let $\Curv_{\pi}^{X}$ be the family of smooth curves of genus $\pi\geq 1$ with an involution over a
smooth curve $X$ of genus $g>0$. Let $r=2(\pi-1)-4(g-1)$.
\begin{enumerate}
\item If $r>0$ then $dim(\Curv_{\pi}^{X})=r-dim(Aut(X))$.

\item If $r=0$ then $dim(\Curv_{\pi}^{X})=0$.

\item If $r<0$ then $\Curv_{\pi}^{X}=\emptyset$.

\end{enumerate}

\end{prop}
{\bf Proof:} Let $\C$ a smooth curve with an involution over the curve $X$ of genus $g$. Let $\gamma:C\lrw
X$ be the
$2:1$ map. By Hurwitz Theorem we know that the map $X$ has $r=2(\pi-1)-4(g-1)$ ramifications. If $r$ is
negative $\Curv_{\pi}^{X}=\emptyset$.

\begin{enumerate}
\item Suppose that $r>0$. Consider the following incidence variety:
$$
J=\{(\C,\R)\in \Curv_{\pi}^{X}\times U^r/ \C \mbox { has an involution branched at
$\R\in X$}\}
$$
where $U_r\subset S^rX$ are the open set of $r$ unordered different points.
We have two projection maps: $p:J\lrw \Curv_{\pi}^{X}$ and $q:J\lrw U^r$. 

Given a curve $\C\in \Curv_{\pi}^{X}$ there is a $2:1$ map $\C\lrw X$ with $r$ ramifications. From
this, $p$ is a surjection. Moreover, since $\pi>1$ the group of automorphisms of $C$ is finite, so there
are a finite number of involutions of $C$ over $X$. In this way we obtain a finite number of
possible ramification points of $X$ up to automorphisms of $X$. But if $g=1$, $dim(Aut(X))=1$ and
$r\geq 1$ and if $g>1$, $dim(Aut(X))=0$. From this there is at most a finite number of automorphism
fixing $r$ generic points. Therefore, $dim(p^{-1}(\C))=dim(Aut(X))$.

On the other hand, given a set $\B\in X$ of $r$ different points we can take a divisor $\b\in Div(X)$
such that $2\b-2\k\sim \B$. Let $S_{\bb}$ be the corresponding canonical ruled surface. By Proposition
\ref{unicidad}, there is a unique curve $\C\in |2X_1|\subset S_{\bb}$ with an involution over $X$ with
branch points over the set
$\B$, so $q$ is a surjection. Furthermore, we know that a curve $\C$ with an involution over the curve
$X$ lays on the linear system
$|2X_1|$ of a canonical ruled surface. Since there are a finite number of divisors $\b$ satisfying
$2\b-2\k\sim \b$ (see \cite{fuentes2} ), we see that $dim(q^{-1}(\C))=0$ and $dim(J)=dim(U^r)=r$.

Thus, we have:
$$
dim(\Curv_{\pi}^{X})=dim(J)-dim(p^{-1}(\C))=r-dim(Aut(X)).
$$

\item Suppose that $r=0$. In this case there are not ramification points. Given a curve $\C\in \Curv_{\pi}^{X}$ we know that lays on the canonical system $|2X_1|$ of a canonical ruled surface $S_{\bb}$
with $2\b-2\k\sim 0$. All curves of this system are isomorphic (Proposition \ref{unicidad}). Moreover,
there are a finite number of divisors $\b$ verifying $2\b\sim 2\k$. Thus $dim(
\Curv_{\pi}^{X})=0$.
 \qed

\end{enumerate}

\begin{prop}\label{dim11}

Let $\CurvHyper_{\pi}^{X}$ be the family of smooth hyperelliptic curves of genus $\pi\geq 1$
with an involution over a smooth curve $X$ of genus $g$. Let $r=2(\pi-1)-4(g-1)$. Then:

\begin{enumerate}

\item If $X$ is neither elliptic nor hyperelliptic then $\CurvHyper_{\pi}^{X}=\emptyset$.

\item If $X$ is elliptic or hyperelliptic then:

\begin{enumerate}

\item If $r>4$ ($\pi>2g+1$) or $r<0$ ($\pi<2g-1$) then $\CurvHyper_{\pi}^{X}=\emptyset$.

\item if $r=4$ ($\pi=2g+1$) then $dim(\CurvHyper_{\pi}^{X})=2$

\item If $r=2$ ($\pi=2g$) then $dim(\CurvHyper_{\pi}^{X})=1$.

\item If $r=0$ ($\pi=2g-1$) then $dim(\CurvHyper_{\pi}^{X})=0$.

\end{enumerate}

\end{enumerate}

\end{prop}

{\bf Proof:} We apply Theorem $3.6$ of \cite{fuentes2}. We see that
$\CurvHyper_{\pi}^{X}=\emptyset$ except when $X$ is elliptic or hyperelliptic and $r=0,2,4$.

\begin{enumerate}
\item Suppose that $r=2$ or $r=4$. By Theorem $3.6$ of \cite{fuentes2}, the branch divisor $\B$ verifies:

\begin{enumerate}

\item If $X$ is hyperelliptic and $r=4$ then $\B\sim 2g^1_2$.

\item If $X$ is elliptic and $r=4$ then $\B\sim a_1+a_2+a_3+a_4$ with $a_1+a_2\sim a_3+a_4$.

\item If $X$ is hyperelliptic and $r=2$ then $\B\sim g^1_2$.

\item If $X$ is elliptic and $r=2$ then $\B\sim a_1+a_2$ for any $a_1,a_2\in X$, $a_1\neq a_2$.

\end{enumerate}

Thus, if we consider the incidence variety:
$$
J_h=\{(\C,\R)\in \CurvHyper_{\pi}^{X}\times U^r/ \C \mbox { has an involution branched at
$\R\in X$}\}
$$
in this case the projection map $q_h:J_h\lrw U^r$ is not a surjection. In fact we have:

\begin{enumerate}

\item If $X$ is hyperelliptic and $r=4$ then $dim(Im(q_h))=2$.

\item If $X$ is elliptic and $r=4$ then $dim(Im(q_h))=3$.

\item If $X$ is hyperelliptic and $r=2$ then $dim(Im(q_h))=1$.

\item If $X$ is elliptic and $r=2$ then $dim(Im(q_h))=2$.

\end{enumerate}

Now, reasoning as in the proposition above we obtain:
$$
dim(\CurvHyper_{\pi}^{X})=dim(J_h)-dim(p^{-1}(\C))=dim(Im(q_h))-dim(Aut(X)).
$$
and the result follows.

\item Suppose that $r=0$. By Theorem $3.6$ of \cite{fuentes2}, $\CurvHyper_{\pi}^{X}\neq \emptyset$ when
$X$ is elliptic or hyperelliptic and then $dim(\CurvHyper_{\pi}^{X})=0$. \qed

\end{enumerate}

\begin{prop}\label{dim2}

Let $\Curv_{\pi}^{g}$ be the family of smooth curves of genus $\pi>1$ with an involution over a curve
 of $g\geq 0$. Let $r=2(\pi-1)-4(g-1)$.
\begin{enumerate}
\item If $r>0$ ($\pi>2g-1$) then $dim(\Curv_{\pi}^{g})=r-dim(Aut(X))+dim({\cal M}_g)$.

\item If $r=0$ ($\pi=2g-1$) then $dim(\Curv_{\pi}^{g})=dim({\cal M}_g)$.

\item If $r<0$ ($\pi<2g-1$) then $dim(\Curv_{\pi}^{g})=\emptyset$.
\end{enumerate}
From this, $dim(\Curv_{\pi}^{g})=2\pi-g-1$. 
Moreover, the dimension of the family of curves with an involution over a hyperelliptic curve of genus
$g$ is $2\pi-2g+1$. 

\end{prop}
{\bf Proof:} Since a curve of genus $\pi>1$ has at most a finite number of involutions, such curve only
has involutions over a finite number of curves of genus $g$. Therefore,
$dim(\Curv_{\pi}^{g})=dim(\Curv_{\pi}^{X})+dim({\cal M}_g)$. 

Moreover, we know that:
$$
\begin{array}{l}
{\mbox{When $g=0$ then $dim({\cal M}_g)=0$ and $dim(Aut(X))=3$.}}\\
{}\\
{\mbox{When $g=1$ then $dim({\cal M}_g)=1$ and $dim(Aut(X))=1$.}}\\
{}\\
{\mbox{When $g>1$ then $dim({\cal M}_g)=3(g-1)$ and $dim(Aut(X))=0$.}}\\
{}\\
{\mbox{If $g\geq 2$ the dimension of the family of hyperelliptic curves of genus $g$ is}}\\
{\mbox{$dim(\Curv_g^0)=2g-1$.}}\\
\end{array}
$$
We have supposed that $\pi>1$ so $r>0$ when $g=1$ or $g=0$.
Now, applying Proposition \ref{dim1} the result follows. \qed

\begin{prop}\label{dim21}

Let $\CurvHyper_{\pi}^{g}$ be the family of smooth hyperelliptic curves of genus $\pi>1$ with an
involution over a curve of $g\geq 0$. Let $r=2(\pi-1)-4(g-1)$.
\begin{enumerate}
\item If $r=0,2,4$ ($\pi=2g-1,2g,2g+1$) then $dim(\CurvHyper_{\pi}^{g})=\pi$.

\item If $r<0$ ($\pi<2g-1$) or $r>4$ ($\pi>2g+1$) then $\Curv_{\pi}^{g}=\emptyset$.
\end{enumerate}

\end{prop}
{\bf Proof:} Since a curve of genus $\pi>1$ has at most a finite number of involutions, such curve only
has involutions over a finite number of curves of genus $g$. Therefore,
$dim(\CurvHyper_{\pi}^{g})=dim(\CurvHyper_{\pi}^{X})+dim(\Curv_{g}^{0})$. Moreover, we know that
$dim(\Curv_{g}^{0})=2g-1$. Applying the Proposition \ref{dim11} the result follows. \qed

\begin{rem}

Note that the dimension of $\CurvHyper_{\pi}^{g}$ does not depend of $g$. The reason is that an
involution of a hyperelliptic curve $C$ of genus $g$ provides an involution of genus $\pi-g$, if we
compose it with the canonical involution. We will study this situation in next section.

\end{rem}

\section{The hyperelliptic case.}

Let $C$ be an hyperelliptic curve of genus $\pi\geq 2$. $C$ has a canonical involution defined by its
unique $\pi_2^1$. We will denote it by $\delta:C\lrw C$, with $\delta(P)=\pi_2^1-P$. Note that any
automorphism $\nu:C\lrw C$ of $C$ verifies that $\nu(\pi_2^1)=\pi_2^1$, so it commutes with $\delta$.

Let $\sigma:C\lrw C$ be an involution of genus $g$ and $\gamma:C\lrw X$ the corresponding double cover.
Let $\b$ the divisor of $X$ such that $\gamma_*\Te_C\sim \Te_X\oplus \Te_X(\k-\b)$.

\begin{lemma}\label{lema1}
Let $P$ be a point of $C$. Then $\delta(P)=\sigma(P)$ if and only if $x=\gamma(\P)$ is a base point of
$\b$.
\end{lemma}
{\bf Proof:}
Let $x=\gamma(\P)$. Then $\delta(P)=\sigma(P)$ if and only if $\gamma^*(x)\sim \pi_2^1$. But
$$
\gamma^*(x)\sim \pi_2^1 \iff h^0(\Te_C(\K_C-\gamma^*(x)))=h^0(\Te_C(\K_C))-1
$$
and
$$
\begin{array}{rl}
{h^0(\Te_C(\K_C-\gamma^*(x)))}&{=h^0(\Te_{S_{\bb}}(X_0+(\b-x)f)=}\\
{}&{=h^0(\Te_X(\K-x))+h^0(\Te_X(\b-x))}\\
{}&{}\\
{h^0(\Te_C(\K_C))}&{=h^0(\Te_{S_{\bb}}(X_0+\b f)=h^0(\Te_X(\K))+h^0(\Te_X(\b))}\\
\end{array}
$$
Because $\K_X$ is base-point-free the conclusion follows. \qed

\begin{teo}
Let $C$ an hyperelliptic curve of genus $\pi\geq 2$. Let $\sigma:C\lrw C$ be an involution of genus $g$.
Then $\delta \sigma$ is an involution of $C$ of genus $\pi-g$.
\end{teo}
{\bf Proof:} Note that $\delta \sigma$ is an involution, because $\sigma$ commutes with any
automorphism of $C$.

Let us study the ramifications points of $\delta \sigma$. We have that $\delta \sigma(P)=P\iff
\sigma(P)=\delta(P)$. By Lemma \ref{lema1}, this happens when $\gamma(x)$ is a base point of $\b$. Thus,
the ramification points of $\delta \sigma$ are $\{\gamma^{-1}(x)/x$ is a base point of $\b\}$.By
Theorem $3.6$ of \cite{fuentes2} we know:

\begin{enumerate}

\item If $\pi=2g+1$, $\b=\k+g_2^1$ and $\b$ is base-point-free.

\item If $\pi=2g$, $\b=\k+P$, with $2P\sim g_2^1$ and $\b$ has one base point.

\item If $\pi=2g-1$, $\b=\sum_1^{g-2} g_2^1+P+Q$, with $2P\sim 2Q\sim g_2^1$, and $\b$ has two base
points except when $X$ is elliptic. But we have supposed $\pi\geq 2$ so in this case $g> 1$.

\end{enumerate}

From this we see that the number of ramifications of $\delta \sigma$ is $0,2,4$ when
$\pi=2g+1,2g,2g-1$ respectively. Applying Hurwitz's formula we obtain that the genus of $\delta
\sigma$ is $\pi-g$. \qed

\begin{cor}
If $\pi\geq 2$ and $g\geq 1$ then $\CurvHyper_{\pi}^g=\CurvHyper_{\pi}^{\pi-g}$. \qed
\end{cor}

\end{document}